\documentclass[11pt]{article}
\usepackage{amssymb}
\usepackage{amsfonts}
\usepackage{amsmath}

\input{tcilatex}
\begin{document}

\begin{center}
\textbf{APPROXIMATE SOLUTIONS FOR MHD SQUEEZING FLUID FLOW BY REPRODUCING
KERNEL HILBERT SPACE METHOD}

MUSTAFA INC, AL\.{I} AKG\"{U}L
\end{center}

\bigskip

\textbf{Abstract:} In this paper, a steady axisymmetric MHD flow of two
dimensional incompressible fluids has been investigated. Reproducing Kernel
Hilbert Space Method (RKHSM) is implemented to obtain solution of reduced
fourth order nonlinear boundary value problem. Numerical results have been
compared with the resutls that obtained by the Range-Kutta Method (RK-4) and
Optimal Homotopy Asymptotic Method (OHAM).

\bigskip

\textbf{AMS Mathematics Subject Classification: }46E22, 35A24

\textbf{Keywords:} Reproducing kernel method, series solutions, squeezing
fluid flow, magnetohydrodynamics, reproducing kernel space.

\bigskip

\section{Introduction}

\ \ \ \ \ \ \ \ Squeezing flows have many applications in food industry,
principally in chemical engineering \cite{papanas}-\cite{ran}. Some
practical examples of squeezing flow include polymer processing, compression
and injection molding. Grimm \cite{grimm} studied numerically, the thin
Newtonian liquids \ films being squeezed between two plates. Squeezing flow
coupled with magnetic field is widely applied to bearing with liquid-metal
lubrication \cite{stefa}, \cite{kamiyama}-\cite{bhattac}. In this paper, we
use RKHSM to study the squeezing MHD fluid flow between two infinite planar
plates.

\ \ \ \ \ Consider a squeezing flow of an incompressible Newtonian fluid in
the presence of a magnetic field of a constant density $\rho $\ and
viscosity $\mu $\ squeezed between two large planar parallel plates,
separated by a small distance $2H$ and the plates approaching each other
with a low constant velocity $V$, as illustrated in Figure 1 and the flow
can be assumed to quasi-steady \cite{papanas}-\cite{ghori}, \cite{idrees}.
The Navier-Stokes equations \cite{ghori}-\cite{ran} governing such flow in
the presence of magnetic field, when inertial terms are retained in the
flow, are given as \cite{islam}%
\begin{equation}
\nabla V.u=0,  \tag{1.1}
\end{equation}%
\begin{equation}
\rho \left[ \frac{\partial u}{\partial t}+\left( u.\nabla \right) u\right]
=\nabla .T+J\times B+\rho f,  \tag{1.2}
\end{equation}%
where $u$ is the velocity vector, $\nabla $ denotes the material time
derivative, $T$ is the Cauchy stress tensor,%
\begin{equation*}
T=-pI+\mu A_{1},
\end{equation*}%
and%
\begin{equation*}
A_{1}=\nabla u+u^{T},
\end{equation*}%
$J$ is the electric current density, $B$ is the total magnetic field and%
\begin{equation*}
B=B_{0}+b,
\end{equation*}%
$B_{0}$ represents the imposed magnetic field and $b$ denotes the induced
magnetic field. In the absence of displacement currents, the modified Ohm's
law and Maxwell's equations (\cite{mohyud} and the references therein) are
given by \cite{islam}%
\begin{equation}
J=\sigma \left[ E+u\times B\right] ,  \tag{1.3}
\end{equation}%
\begin{equation}
\func{div}B=0,\text{ \ \ }\nabla \times B=\mu _{m}J,\text{ \ \ }\func{curl}E=%
\frac{\partial B}{\partial t},  \tag{1.4}
\end{equation}%
in which $\sigma $\ is the electrical conductivity, $E$ the electric field
and $\mu _{m}$ the magnetic permeability.

We need the following assumptions \cite{islam}:

$a)$ The density $\rho ,$ magnetic permeability $\mu _{m}$ and electric
field conductivity $\sigma $, are assumed to be constant throughout the flow
field region.

$b)$ The electrical conductivity $\sigma $ of the fluid considers being
finite.

$c)$ Total magnetic field $B$ is perpendicular to the velocity field $V$ and
the induced magnetic field $b$ is negligible compared with the applied
magnetic field $B_{0}$ so that the magnetic Reynolds number is small (\cite%
{mohyud} and the references therein).

$d)$ We assume that a situation where no energy is added or extracted from
the fluid by the electric field, which implies that there is no electric
field present in the fluid flow region.

\bigskip

\begin{equation*}
\FRAME{itbpF}{4.6311in}{2.5339in}{0in}{}{}{Figure}{\special{language
"Scientific Word";type "GRAPHIC";display "USEDEF";valid_file "T";width
4.6311in;height 2.5339in;depth 0in;original-width 7.043in;original-height
14.5453in;cropleft "0";croptop "1";cropright "1";cropbottom "0";tempfilename
'MRCOGT03.wmf';tempfile-properties "XPR";}}
\end{equation*}

\begin{center}
Figure 1. A steady squeezing axisymmetric fluid flow between two parallel
plates.
\end{center}

Under these assumptions, the magnetohydrodynamic force involved in Eq. (1.2)
can be put into the form%
\begin{equation}
J\times B=-\sigma B_{0}^{2}u.  \tag{1.5}
\end{equation}

We consider an incompressible Newtonian fluid, squeezed between two large
planar, parallel smooth plates which is separated by a small distance $2H$
and moving towards each other with velocity $V$. We assume that the plates
are nonconducting and the magnetic field is applied along the z-axis. For
small values of the velocity $V$, as shown in the Figure 1, the gap distance 
$2H$ between the plates changes slowly with time $t$, so that it can be
taken as constant, the flow is steady \cite{stefa},\cite{idrees}. An
axisymmetric flow in cylindrical coordinates $r,\theta ,z$ with $z$-axis
perpendicular to plates and $z=\pm H$ at the plates. Since we have axial
symmetry, so $u$ is represented by%
\begin{equation*}
u=\left( u_{r}\left( r,z\right) ,0,u_{z}\left( r,z\right) \right) ,
\end{equation*}%
when body forcesm are negligible, Navier-Stokes Eqs. (1.1)-(1.2) in
cylindrical coordinates, where there is no tangential velocity $(u_{\theta
}=0)$, are given as \cite{islam}%
\begin{equation}
\rho \left( u_{r}\frac{\partial u_{r}}{\partial r}+u_{z}\frac{\partial u_{r}%
}{\partial z}\right) =-\frac{\partial p}{\partial r}+\left( \frac{\partial
^{2}u_{r}}{\partial r^{2}}+\frac{1}{r}\frac{\partial u_{r}}{\partial r}-%
\frac{u_{r}}{r^{2}}+\frac{\partial ^{2}u_{r}}{\partial z^{2}}\right) +\sigma
B_{0}^{2}u_{r,}  \tag{1.6}
\end{equation}%
\begin{equation}
\rho \left( u_{z}\frac{\partial u_{z}}{\partial r}+u_{z}\frac{\partial u_{z}%
}{\partial z}\right) =-\frac{\partial p}{\partial r}+\left( \frac{\partial
^{2}u_{z}}{\partial r^{2}}+\frac{1}{r}\frac{\partial u_{z}}{\partial r}+%
\frac{\partial ^{2}u_{z}}{\partial z^{2}}\right) ,  \tag{1.7}
\end{equation}%
where $p$ is the pressure, and equation of continuity is:%
\begin{equation}
\frac{1}{r}\frac{\partial }{\partial r}(ru_{r})+\frac{\partial u_{z}}{%
\partial z}=0.  \tag{1.8}
\end{equation}%
The boundary conditions require%
\begin{equation}
\left. 
\begin{array}{c}
u_{r}=0,\text{ \ \ }u_{z}=-V\text{ \ \ at \ \ }z=H, \\ 
\\ 
\frac{\partial u_{r}}{\partial z}=0,\text{ \ \ }u_{z}=0\text{ \ \ at \ \ }%
z=0.%
\end{array}%
\right.  \tag{1.9}
\end{equation}%
Introducing the axisymmetric Stokes stream function $\Psi $ as%
\begin{equation}
u_{r}=\frac{1}{r}\frac{\partial \Psi }{\partial z},\text{ \ \ }u_{z}=-\frac{1%
}{r}\frac{\partial \Psi }{\partial r}.  \tag{1.10}
\end{equation}%
The continuity equation is satisfied using Eq. (1.10). Substituting Eqs.
(1.3)-(1.5) and Eq. (1.10) into the Eqs. (1.7)-(1.8), we obtain%
\begin{equation}
-\frac{\rho }{r^{2}}\frac{\partial \Psi }{\partial r}E^{2}\Psi =-\frac{%
\partial p}{\partial r}+\frac{\mu }{r}\frac{\partial E^{2}\Psi }{\partial z}-%
\frac{\sigma B_{0}^{2}}{r}\frac{\partial \Psi }{\partial z}  \tag{1.11}
\end{equation}%
and%
\begin{equation}
-\frac{\rho }{r^{2}}\frac{\partial \Psi }{\partial z}E^{2}\Psi =-\frac{%
\partial p}{\partial z}+\frac{\mu }{r}\frac{\partial E^{2}\Psi }{\partial r}.
\tag{1.12}
\end{equation}%
Eliminating the pressure from Eqs. (1.11) and (1.12) by integribility
condition we get the compatibility equation as \cite{islam}%
\begin{equation}
-\rho \left[ \frac{\partial \left( \Psi ,\frac{E^{2}\Psi }{r^{2}}\right) }{%
\partial (r,z)}\right] =\frac{\mu }{r}E^{2}\Psi -\frac{\sigma B_{0}^{2}}{r}%
\frac{\partial ^{2\Psi }}{\partial z^{2}},  \tag{1.13}
\end{equation}%
where%
\begin{equation*}
E^{2}=\frac{\partial ^{2}}{\partial r^{2}}-\frac{1}{r}\frac{\partial }{%
\partial r}+\frac{\partial ^{2}}{\partial z^{2}}.
\end{equation*}%
The stream function can be expressed as \cite{papanas}, \cite{ghori}%
\begin{equation}
\Psi (r,z)=r^{2}F(z).  \tag{1.14}
\end{equation}%
In view of Eq. (1.14), the compatibility Eq. (1.13) and the boundary
conditions (1.9) take the form%
\begin{equation}
F^{(iv)}(z)-\frac{\sigma B_{0}^{2}}{r}F^{\prime \prime }(z)+2\frac{\rho }{%
\mu }F(z)F^{\prime \prime \prime }(z)=0,  \tag{1.15}
\end{equation}%
subject to%
\begin{equation}
\left. 
\begin{array}{c}
F(0)=0,\text{ \ \ }F^{\prime \prime }(0)=0, \\ 
\\ 
F(H)=\frac{V}{2},\text{ \ \ }F^{\prime }(H)=0.%
\end{array}%
\right.  \tag{1.16}
\end{equation}%
Introducing the following non-dimensional parameters%
\begin{equation*}
F^{\ast }=2\frac{F}{V},\text{ \ \ }z^{\ast }=\frac{z}{H},\text{ \ \ }\func{Re%
}=\frac{\rho HV}{\mu },\text{ \ \ }m=B_{0}H\sqrt{\frac{\sigma }{\mu }}.
\end{equation*}%
For simplicity omitting the $\ast $, the boundary value problem
(1.15)-(1.16) becomes \cite{islam}%
\begin{equation}
F^{(iv)}(z)-m^{2}F^{\prime \prime }(z)+\func{Re}F(z)F^{\prime \prime \prime
}(z)=0,  \tag{1.17}
\end{equation}%
with the boundary conditions%
\begin{equation}
\left. 
\begin{array}{c}
F(0)=0,\text{ \ \ }F^{\prime \prime }(0)=0, \\ 
\\ 
F(1)=1,\text{ \ \ }F^{\prime }(1)=0,%
\end{array}%
\right.  \tag{1.18}
\end{equation}%
where $\func{Re}$ is the Reynolds number and $m$ is Hartmann number. This
problem has been solved by RKHSM and for comparison it has been compared
with the OHAM and numerically with the RK-4 by using maple 16.

The RKHSM which accurately computes the series solution is of great interest
to applied sciences. The method provides the solution in a rapidly
convergent series with components that can be elegantly computed. The
efficiency of the method was used by many authors to investigate several
scientific applications. Geng and Cui \cite{geng} applied the RKHSM to
handle the second-order boundary value problems. Yao and Cui \cite{yao} and
Wang et al. \cite{wang} investigated a class of singular boundary value
problems by this method and the obtained results were good. Zhou et al. \cite%
{zhou} used the RKHSM effectively to solve second-order boundary value
problems. Wang and Chao \cite{wang1} Li and Cui \cite{li}, Zhou and Cui \cite%
{zhou1} independently employed the RKSHSM to variable-coefficient partial
differential equations. Geng and Cui \cite{geng1}, Du and Cui \cite{du}
investigated the approximate solution of the forced Duffing equation with
integral boundary conditions by combining the homotopy perturbation method
and the RKM. Lv and Cui [32] presented a new algorithm to solve linear
fifth-order boundary value problems. Cui and Du \cite{du1} obtained the
representation of the exact solution for the nonlinear Volterra-Fredholm
integral equations by using the reproducing kernel Hilbert space method. Wu
and Li \cite{wu} applied iterative reproducing kernel Hilbert space method
to obtain the analytical approximate solution of a nonlinear oscillator with
discontinuities. For more details about RKHSM and the modified forms and its
effectiveness, see \cite{cui}-\cite{ali4} and the references therein.\newline

The paper is organized as follows. Section 2 introduces several reproducing
kernel spaces and a linear operator. S{\normalsize olution representation in}%
\textbf{\ }$W_{2}^{5}[0,1]$ has been presented in Section 3. It provides the
main results, the exact and approximate solution of $(1.1)$ and an iterative
method are developed for the kind of problems in the reproducing kernel
space. We have proved that the approximate solution converges to the exact
solution uniformly. Some numerical experiments are illustrated in Section 4.
There are some conclusions in the last section.

\bigskip

\section{Preliminaries}

\textbf{2.1. Reproducing Kernel Spaces}

\bigskip

In this section, we define some useful reproducing kernel spaces.

\textbf{Definition 2.1.} \textit{(Reproducing kernel)}. Let $E$ be a
nonempty abstract set. A function $K:E\times E\longrightarrow C$ is a
reproducing kernel of the Hilbert space $H$ if and only if%
\begin{equation}
\left\{ 
\begin{array}{c}
\forall t\in E,\text{ }K\left( .,t\right) \in H, \\ 
\forall t\in E,\text{ }\forall \varphi \in H,\text{ }\left\langle \varphi
\left( .\right) ,K\left( .,t\right) \right\rangle =\varphi \left( t\right) .%
\end{array}%
\right.  \tag{2.1}
\end{equation}

The last condition is called "the reproducing property": the value of the
function $\varphi $ at the point $t$ is reproduced by the inner product of $%
\varphi $ with $K\left( .,t\right) $

\bigskip

\textbf{Definition 2.2. }We define the space $W_{2}^{5}[0,1]$ by 
\begin{equation*}
W_{2}^{5}[0,1]=\left\{ 
\begin{array}{c}
u\Bigl|u,\text{ }u^{\prime },\text{ }u^{\prime \prime },u^{\prime \prime
\prime },u^{(4)}\text{ are absolutely continuous in }[0,1], \\ 
\\ 
u^{\left( 5\right) }\in L^{2}[0,1],\text{ }x\in \lbrack 0,1],\text{ }%
u(0)=u(1)=u^{\prime }(1)=u^{\prime \prime }(0)=0.%
\end{array}%
\right\}
\end{equation*}%
The fifth derivative of $u$ exists almost everywhere since $u^{(4)}$ is
absolutely continuous. The inner product and the norm in $W_{2}^{5}[0,1]$
are defined respectively by%
\begin{equation*}
\left\langle u,v\right\rangle _{{\LARGE W}_{2}^{5}}=%
\sum_{i=0}^{4}u^{(i)}(0)v^{(i)}(0)+\int_{0}^{1}u^{(5)}(x)v^{(5)}(x)dx,\text{
\ \ }u,v\in W_{2}^{5}[0,1],
\end{equation*}%
and

\begin{equation*}
\left\Vert u\right\Vert _{W_{2}^{5}}=\sqrt{\left\langle u,u\right\rangle
_{_{W_{2}^{5}}}},\ \text{\ \ }u\in W_{2}^{5}[0,1].
\end{equation*}%
The space $W_{2}^{5}[0,1]$ \ is a reproducing kernel space, i.e., for each
fixed $y\in \lbrack 0,1]$ \ and any $u\in W_{2}^{5}[0,1],$ there exists a
function $R_{y}$ such that

\begin{equation*}
u=\left\langle u,R_{y}\right\rangle _{W_{2}^{5}}.
\end{equation*}

\bigskip

\textbf{Definition 2.3.}\ We define the space $W_{2}^{4}[0,1]$ by 
\begin{equation*}
W_{2}^{4}[0,1]=\left\{ 
\begin{array}{c}
u\Bigl|u,\text{ }u^{\prime },\text{ }u^{\prime \prime },\text{ }u^{\prime
\prime \prime }\text{ are absolutely continuous in }[0,1], \\ 
\\ 
u^{\left( 4\right) }\in L^{2}[0,1],\text{ \ \ }x\in \lbrack 0,1].%
\end{array}%
\right\}
\end{equation*}%
The fourth derivative of $u$ exists almost everywhere since $u^{\prime
\prime \prime }$ is absolutely continuous. The inner product and the norm in 
$W_{2}^{4}[0,1]$ are defined respectively by%
\begin{equation*}
\left\langle u,v\right\rangle _{{\LARGE W}_{2}^{4}}=%
\sum_{i=0}^{3}u^{(i)}(0)v^{(i)}(0)+\int_{0}^{1}u^{(4)}(x)v^{(4)}(x)dx,\text{
\ \ }u,v\in W_{2}^{4}[0,1],
\end{equation*}%
and

\begin{equation*}
\left\Vert u\right\Vert _{W_{2}^{4}}=\sqrt{\left\langle u,u\right\rangle
_{_{W_{2}^{4}}}},\ \text{\ \ }u\in W_{2}^{4}[0,1].
\end{equation*}%
The space $W_{2}^{4}[0,1]$ \ is a reproducing kernel space, i.e., for each
fixed $y\in \lbrack 0,1]$ \ and any $u\in W_{2}^{4}[0,1],$ there exists a
function $r_{y}$ such that

\begin{equation*}
u=\left\langle u,r_{y}\right\rangle _{W_{2}^{4}}.
\end{equation*}

\bigskip

\textbf{Theorem 2.1. }The space $W_{2}^{5}[0,1]$ is a reproducing kernel
Hilbert space whose reproducing kernel function is given by%
\begin{equation*}
R_{y}\left( x\right) =\left\{ 
\begin{array}{c}
\displaystyle{\sum_{i=1}^{10}c_{i}\left( y\right) x^{i-1}},\text{ \ \ }x\leq
y, \\ 
\\ 
\displaystyle{\sum_{i=1}^{10}d_{i}\left( y\right) x^{i-1}},\text{ \ \ }x>y,%
\end{array}%
\right.
\end{equation*}%
where, $c_{j}(y)$ can be deduced easily by using for example MAPLE 16, 
\begin{equation*}
c_{1}(y)=0,\text{ \ \ }c_{3}(y)=0,
\end{equation*}
\begin{eqnarray*}
c_{2}(y) &=&-\frac{1537}{378141715}y^{9}+\frac{9374}{378141715}y^{8}-\frac{%
3932}{75628343}y^{7}+\frac{608}{54020245}y^{6} \\
&& \\
&&+\frac{8006}{54020245}y^{5}+\frac{8006}{10804049}y^{4}-\frac{14592}{%
10804049}y^{3}+\frac{5201}{10804049}y,
\end{eqnarray*}

\begin{eqnarray*}
c_{4}(y) &=&\frac{16061}{13613101740}y^{9}+\frac{38}{1134425145}y^{8}-\frac{%
48487}{1134425145}y^{7}+\frac{10758}{54020245}y^{6} \\
&& \\
&&-\frac{145157}{324121470}y^{5}-\frac{145157}{64824294}y^{4}+\frac{1509137}{%
388945764}y^{3}-\frac{14592}{10804049}y,
\end{eqnarray*}%
\begin{eqnarray*}
c_{5}(y) &=&\frac{1243}{10890481392}y^{9}-\frac{4003}{217809627840}y^{8}-%
\frac{107867}{27226203480}y^{7}+\frac{145157}{7778915280}y^{6} \\
&& \\
&&-\frac{81901}{1944728820}y^{5}+\frac{9493633}{6223132224}y^{4}-\frac{145157%
}{64824294}y^{3}+\frac{8006}{10804049}y,
\end{eqnarray*}%
\begin{eqnarray*}
c_{6}(y) &=&\frac{1243}{54452406960}y^{9}-\frac{4003}{1089048139200}y^{8}-%
\frac{107867}{136131017400}y^{7}+\frac{145157}{38894576400}y^{6} \\
&& \\
&&-\frac{81901}{9723644100}y^{5}+\frac{9493633}{31115661120}y^{4}-\frac{%
145157}{324121470}y^{3}+\frac{8006}{54020245}y,
\end{eqnarray*}%
\begin{eqnarray*}
c_{7}(y) &=&-\frac{16061}{1633572208800}y^{9}-\frac{19}{68065508700}y^{8}+%
\frac{48487}{136131017400}y^{7}-\frac{1793}{1080404900}y^{6} \\
&& \\
&&+\frac{145157}{38894576400}y^{5}+\frac{145157}{7778915280}y^{4}-\frac{%
1509137}{46673491680}y^{3}+\frac{608}{54020245}y,
\end{eqnarray*}

\begin{eqnarray*}
c_{8}(y) &=&\frac{3631}{1905834243600}y^{9}+\frac{983}{762333697440}y^{8}-%
\frac{18797}{238229280450}y^{7}+\frac{48487}{136131017400}y^{6} \\
&& \\
&&-\frac{107867}{136131017400}y^{5}-\frac{107867}{27226203480}y^{4}-\frac{%
48487}{1134425145}y^{3}+\frac{1}{10080}y^{2}-\frac{3932}{75628343}y,
\end{eqnarray*}%
\begin{eqnarray*}
c_{9}(y) &=&\frac{1537}{15246673948800}y^{9}-\frac{4687}{7623336974400}y^{8}+%
\frac{983}{7623336974400}y^{7}-\frac{19}{68065508700}y^{6} \\
&& \\
&&-\frac{4003}{1089048139200}y^{5}-\frac{4003}{217809627840}y^{4}+\frac{38}{%
1134425145}y^{3}-\frac{743}{62231322240}y,
\end{eqnarray*}

\begin{eqnarray*}
c_{10}(y) &=&-\frac{323}{4288127048100}y^{9}+\frac{1537}{15246673948800}%
y^{8}+\frac{3631}{1905834243600}y^{7}-\frac{16061}{1633572208800}y^{6} \\
&& \\
&&+\frac{1243}{54452406960}y^{5}+\frac{1243}{10890481392}y^{4}+\frac{16061}{%
13613101740}y^{3}-\frac{1537}{378141715}y+\frac{1}{362880},
\end{eqnarray*}%
Similarly, we can obtain, 
\begin{equation*}
d_{1}(y)=\frac{1}{362880}y^{9},\text{ \ \ \ }d_{3}(y)=\frac{1}{10080}y^{7},
\end{equation*}%
\begin{eqnarray*}
d_{2}(y) &=&-\frac{1537}{378141715}y^{9}-\frac{743}{62231322240}y^{8}-\frac{%
3932}{75628343}y^{7}+\frac{608}{54020245}y^{6} \\
&& \\
&&+\frac{8006}{54020245}y^{5}+\frac{8006}{10804049}y^{4}-\frac{14592}{%
10804049}y^{3}+\frac{5201}{10804049}y,
\end{eqnarray*}

\begin{eqnarray*}
d_{4}(y) &=&\frac{16061}{13613101740}y^{9}+\frac{38}{1134425145}y^{8}-\frac{%
48487}{1134425145}y^{7}-\frac{1509137}{46673491680}y^{6} \\
&& \\
&&-\frac{145157}{324121470}y^{5}-\frac{145157}{64824294}y^{4}+\frac{1509137}{%
388945764}y^{3}-\frac{14592}{10804049}y,
\end{eqnarray*}%
\begin{eqnarray*}
d_{5}(y) &=&\frac{1243}{10890481392}y^{9}-\frac{4003}{217809627840}y^{8}-%
\frac{107867}{27226203480}y^{7}+\frac{145157}{7778915280}y^{6} \\
&& \\
&&+\frac{9493633}{31115661120}y^{5}+\frac{9493633}{6223132224}y^{4}-\frac{%
145157}{64824294}y^{3}+\frac{8006}{10804049}y,
\end{eqnarray*}%
\begin{eqnarray*}
d_{6}(y) &=&\frac{1243}{54452406960}y^{9}-\frac{4003}{1089048139200}y^{8}-%
\frac{107867}{136131017400}y^{7}+\frac{145157}{38894576400}y^{6} \\
&& \\
&&-\frac{81901}{9723644100}y^{5}-\frac{81901}{1944728820}y^{4}-\frac{145157}{%
324121470}y^{3}+\frac{8006}{54020245}y,
\end{eqnarray*}%
\begin{eqnarray*}
d_{7}(y) &=&-\frac{16061}{1633572208800}y^{9}-\frac{19}{68065508700}y^{8}+%
\frac{48487}{136131017400}y^{7}-\frac{1793}{1080404900}y^{6} \\
&& \\
&&+\frac{145157}{38894576400}y^{5}+\frac{145157}{7778915280}y^{4}+\frac{10758%
}{54020245}y^{3}+\frac{608}{54020245}y,
\end{eqnarray*}

\begin{eqnarray*}
d_{8}(y) &=&\frac{3631}{1905834243600}y^{9}+\frac{983}{762333697440}y^{8}-%
\frac{18797}{238229280450}y^{7}+\frac{48487}{136131017400}y^{6} \\
&& \\
&&-\frac{107867}{136131017400}y^{5}-\frac{107867}{27226203480}y^{4}-\frac{%
48487}{1134425145}y^{3}-\frac{3932}{75628343}y,
\end{eqnarray*}%
\begin{eqnarray*}
d_{9}(y) &=&\frac{1537}{15246673948800}y^{9}-\frac{4687}{7623336974400}y^{8}+%
\frac{983}{7623336974400}y^{7}-\frac{19}{68065508700}y^{6} \\
&& \\
&&-\frac{4003}{1089048139200}y^{5}-\frac{4003}{217809627840}y^{4}+\frac{38}{%
1134425145}y^{3}+\frac{9374}{378141715}y,
\end{eqnarray*}

\begin{eqnarray*}
d_{10}(y) &=&-\frac{323}{4288127048100}y^{9}+\frac{1537}{15246673948800}%
y^{8}+\frac{3631}{1905834243600}y^{7}-\frac{16061}{1633572208800}y^{6} \\
&& \\
&&+\frac{1243}{54452406960}y^{5}+\frac{1243}{10890481392}y^{4}+\frac{16061}{%
13613101740}y^{3}-\frac{1537}{378141715}y.
\end{eqnarray*}

\bigskip

\textbf{Proof: }Let $u\in W_{2}^{5}[0,1].$ By Definition 2.2 we have 
\begin{equation}
\left\langle u,{\normalsize R}_{y}\right\rangle _{_{{\LARGE W}%
_{2}^{5}}}=\sum_{i=0}^{4}u^{(i)}(0){\normalsize R}_{y}^{(i)}{\normalsize (0)}%
+\int_{0}^{1}u^{(5)}(x){\normalsize R}_{y}^{(5)}{\normalsize (x)}dx. 
\tag{2.3}
\end{equation}%
Through several integrations by parts for (2.3) we have

\begin{eqnarray}
\left\langle u,{\normalsize R}_{y}\right\rangle _{_{W_{2}^{5}}}
&=&\sum_{i=0}^{4}{\normalsize u}^{(i)}{\normalsize (0)}\left[ {\normalsize R}%
_{y}^{(i)}{\normalsize (0)-(-1)}^{(4-i)}{\normalsize R}_{y}^{(9-i)}%
{\normalsize (0)}\right]  \TCItag{2.4} \\
&&{\normalsize +}\sum_{i=0}^{4}{\normalsize (-1)}^{(4-i)}{\normalsize u}%
^{(i)}{\normalsize (1)R}_{y}^{(9-i)}{\normalsize (1)-}\int_{0}^{1}u%
{\normalsize (x)R}_{y}^{(10)}{\normalsize (x)dx.}  \notag
\end{eqnarray}%
Note that property of the reproducing kernel

\begin{equation*}
\left\langle u,R_{y}\right\rangle _{W_{2}^{5}}=u(y).
\end{equation*}%
Now, if 
\begin{equation}
\left\{ 
\begin{array}{c}
R_{y}^{^{\prime }}(0)+R_{y}^{(8)}(0)=0, \\ 
R_{y}^{(3)}(0)+R_{y}^{(6)}(0)=0, \\ 
R_{y}^{(4)}(0)-R_{y}^{(5)}(0)=0, \\ 
R_{y}^{(5)}(1)=0, \\ 
R_{y}^{(6)}(1)=0, \\ 
R_{y}^{(7)}(1)=0,%
\end{array}%
\right.  \tag{2.5}
\end{equation}%
then (2.4) implies that,

\begin{equation*}
{\normalsize R}_{y}^{(10)}{\normalsize (x)=-\delta (x-y),}
\end{equation*}%
when $x\neq y,$ then

\begin{equation*}
{\normalsize R}_{y}^{(10)}{\normalsize (x)=0,}
\end{equation*}%
and therefore

\begin{equation*}
{\normalsize R}_{y}{\normalsize (x)=}\left\{ 
\begin{array}{c}
\displaystyle{\ \sum_{i=1}^{10}c_{i}(y)x^{i-1}},\text{ \ \ }x\leq y, \\ 
\\ 
\displaystyle{\sum_{i=1}^{10}d_{i}(y)x^{i-1}},\text{ \ \ }x>y.%
\end{array}%
\right.
\end{equation*}%
Since

\begin{equation*}
{\normalsize R}_{y}^{(10)}{\normalsize (x)=\delta (x-y),}
\end{equation*}%
we have

\begin{equation}
{\normalsize \partial }^{k}{\normalsize R}_{y^{+}}{\normalsize (y)=\partial }%
^{k}{\normalsize R}_{y^{-}}{\normalsize (y),}\text{ \ }k=0,1,2,3,4,5,6,7,8, 
\tag{2.6}
\end{equation}%
and

\begin{equation}
{\normalsize \partial }^{9}{\normalsize R}_{y^{{\LARGE +}}}{\normalsize %
(y)-\partial }^{9}{\normalsize R}_{y^{-}}{\normalsize (y)=-1.}  \tag{2.7}
\end{equation}%
Since ${\normalsize R}_{y}{\normalsize (x)\in }W_{2}^{5}[0,1]$, it follows
that

\begin{equation}
R_{y}(0)=0,\text{ }R_{y}(1)=0,\text{ }R_{y}^{\prime }(1)=0,\text{ }%
R_{y}^{\prime \prime }(0)=0\text{.}  \tag{2.8}
\end{equation}%
From (2.5)-(2.8), the unknown coefficients $c_{i}(y)$ and $d_{i}(y)$ $%
(i=1,2,...,12)$ can be obtained. This completes the proof. \ \ \ \ \ \ $%
\square $

\textbf{\bigskip }

\textbf{Remark 2.1. }\ Reprroducing kernel function $r_{y}$ of $%
W_{2}^{4}[0,1]$ is given as%
\begin{equation*}
r_{y}(x)=\left\{ 
\begin{array}{c}
1+xy+\frac{1}{4}y^{2}x^{2}+\frac{1}{36}y^{3}x^{3}+\frac{1}{144}y^{3}x^{4}-%
\frac{1}{240}y^{2}x^{5}+\frac{1}{720}yx^{6}-\frac{1}{5040}x^{7},\text{ \ \ }%
x\leq y, \\ 
\\ 
1+yx+\frac{1}{4}y^{2}x^{2}+\frac{1}{36}y^{3}x^{3}+\frac{1}{144}x^{3}y^{4}-%
\frac{1}{240}x^{2}y^{5}+\frac{1}{720}xy^{6}-\frac{1}{5040}y^{7},\text{ \ \ }%
x>y.%
\end{array}%
\right.
\end{equation*}%
This can be proved easily as in Theorem 2.1.

\bigskip

\section{Solution Representation in $W_{2}^{5}[0,1]$}

In this section, the solution of equation (1.17) is given in the reproducing
kernel space $W_{2}^{5}[0,1].$

On defining the linear operator $L:W_{2}^{5}[0,1]\rightarrow W_{2}^{4}[0,1]$%
\ as

\begin{equation*}
Lu=u^{(4)}(x)+\func{Re}\frac{e^{x}}{e}\left( x^{3}-4x^{2}+4x\right)
-m^{2}u^{\prime \prime }(x)+\func{Re}\frac{e^{x}}{e}\left(
x^{3}+5x^{2}-2x-6\right) u(x).
\end{equation*}%
Model problem (1.17) changes the following problem:

\begin{equation}
\left\{ 
\begin{array}{c}
Lu=M(x,u,u^{(3)}),,\text{ \ \ }x\in \lbrack 0,1], \\ 
\\ 
u(0)=0,\text{ \ \ }u(1)=0,\text{ \ \ \ }u^{\prime }(1)=0,\text{ \ \ }%
u^{\prime \prime }(0)=0,%
\end{array}%
\right.  \tag{3.1}
\end{equation}%
where%
\begin{eqnarray*}
M(x,u,u^{(3)}) &=&-\func{Re}u^{(3)}(x)u(x)-\func{Re}\left( \frac{e^{x}}{e}%
\right) ^{2}\left( x^{3}-4x^{2}+4x\right) \left( x^{3}+5x^{2}-2x-6\right) \\
&&-\frac{e^{x}}{e}\left( x^{3}+8x^{2}+8x-2\right) +m^{2}\frac{e^{x}}{e}%
\left( x^{3}+2x^{2}-6x\right) .
\end{eqnarray*}

\subsection{\textbf{The Linear Boundedness of Operator }$L.$}

\bigskip

\textbf{Theorem 3.1. }The operator $L$ defined by (3.1) is a bounded linear
operator.

\textbf{Proof: }We only need to prove%
\begin{equation*}
\left\Vert Lu\right\Vert _{{\LARGE W}_{{\LARGE 2}}^{{\LARGE 4}}}^{2}\leq
P\left\Vert Lu\right\Vert _{{\LARGE W}_{{\LARGE 2}}^{{\LARGE 5}}}^{2},
\end{equation*}%
where $P$ is a positive constant. By Definition 2.3, we have%
\begin{equation*}
\left\Vert u\right\Vert _{{\LARGE W}_{2}^{4}}^{2}=\left\langle
u,u\right\rangle _{{\LARGE W}_{2}^{4}}=\sum_{i=0}^{3}\left[ u^{(i)}(0)\right]
^{2}+\int_{0}^{1}\left[ u^{(4)}(x)\right] ^{2}dx,\text{ \ \ }u\in
W_{2}^{4}[0,1],
\end{equation*}%
and

\begin{eqnarray*}
\left\Vert Lu\right\Vert _{{\LARGE W}_{{\LARGE 2}}^{{\LARGE 4}}}^{2}
&=&\left\langle Lu,Lu\right\rangle _{{\LARGE W}_{{\LARGE 2}}^{{\LARGE 4}}}=%
\left[ \left( Lu\right) \left( 0\right) \right] ^{2}+\left[ \left( Lu\right)
^{\prime }\left( 0\right) \right] ^{2}+\left[ \left( Lu\right) ^{\prime
\prime }\left( 0\right) \right] ^{2} \\
&&+\left[ \left( Lu\right) ^{(3)}\left( 0\right) \right] ^{2}+\int_{0}^{1}%
\left[ \left( Lu\right) ^{(4)}\left( x\right) \right] ^{2}dx.
\end{eqnarray*}%
By reproducing property, we have%
\begin{equation*}
u\left( x\right) =\left\langle u,R_{x}\right\rangle _{{\LARGE W}_{{\LARGE 2}%
}^{{\LARGE 5}}},
\end{equation*}%
and%
\begin{eqnarray*}
\left( Lu\right) \left( x\right) &=&\left\langle u,\left( LR_{x}\right)
\right\rangle _{{\LARGE W}_{{\LARGE 2}}^{{\LARGE 5}}},\text{ \ \ }\left(
Lu\right) ^{\prime }\left( x\right) =\left\langle u,\left( LR_{x}\right)
^{\prime }\right\rangle _{{\LARGE W}_{{\LARGE 2}}^{{\LARGE 5}}}, \\
&& \\
\left( Lu\right) ^{\prime \prime }\left( x\right) &=&\left\langle u,\left(
LR_{x}\right) ^{\prime \prime }\right\rangle _{{\LARGE W}_{{\LARGE 2}}^{%
{\LARGE 5}}},\text{ \ \ }\left( Lu\right) ^{(3)}\left( x\right)
=\left\langle u,\left( LR_{x}\right) ^{(3)}\right\rangle _{{\LARGE W}_{%
{\LARGE 2}}^{{\LARGE 5}}}, \\
&& \\
\left( Lu\right) ^{(4)}\left( x\right) &=&\left\langle u,\left(
LR_{x}\right) ^{(4)}\right\rangle _{{\LARGE W}_{{\LARGE 2}}^{{\LARGE 5}}}.
\end{eqnarray*}%
Therefore%
\begin{eqnarray*}
\left\vert \left( Lu\right) \left( x\right) \right\vert &\leq &\left\Vert
u\right\Vert _{{\LARGE W}_{{\LARGE 2}}^{{\LARGE 5}}}\left\Vert
LR_{x}\right\Vert _{W_{{\LARGE 2}}^{{\LARGE 5}}}\leq a_{1}\left\Vert
u\right\Vert _{{\LARGE W}_{{\LARGE 2}}^{{\LARGE 5}}},\text{ \ (where }a_{1}>0%
\text{ is a positive constant),} \\
&& \\
\left\vert \left( Lu\right) ^{\prime }\left( x\right) \right\vert &\leq
&\left\Vert u\right\Vert _{{\LARGE W}_{{\LARGE 2}}^{{\LARGE 5}}}\left\Vert
\left( LR_{x}\right) ^{\prime }\right\Vert _{W_{{\LARGE 2}}^{{\LARGE 5}%
}}\leq a_{2}\left\Vert u\right\Vert _{{\LARGE W}_{{\LARGE 2}}^{{\LARGE 5}}},%
\text{ \ (where }a_{2}>0\text{ is a positive constant),} \\
&& \\
\left\vert \left( Lu\right) ^{\prime \prime }\left( x\right) \right\vert
&\leq &\left\Vert u\right\Vert _{{\LARGE W}_{{\LARGE 2}}^{{\LARGE 5}%
}}\left\Vert \left( LR_{x}\right) ^{\prime \prime }\right\Vert _{W_{{\LARGE 2%
}}^{{\LARGE 5}}}\leq a_{3}\left\Vert u\right\Vert _{{\LARGE W}_{{\LARGE 2}}^{%
{\LARGE 5}}},\text{ \ (where }a_{3}>0\text{ is a positive constant),} \\
&& \\
\left\vert \left( Lu\right) ^{(3)}\left( x\right) \right\vert &\leq
&\left\Vert u\right\Vert _{{\LARGE W}_{{\LARGE 2}}^{{\LARGE 5}}}\left\Vert
\left( LR_{x}\right) ^{(3)}\right\Vert _{W_{{\LARGE 2}}^{{\LARGE 5}}}\leq
a_{4}\left\Vert u\right\Vert _{{\LARGE W}_{{\LARGE 2}}^{{\LARGE 5}}},\text{
\ (where }a_{4}>0\text{ is a positive constant),}
\end{eqnarray*}%
Thus

\begin{equation*}
\left[ \left( Lu\right) \left( 0\right) \right] ^{2}+\left[ \left( Lu\right)
^{\prime }\left( 0\right) \right] ^{2}+\left[ \left( Lu\right) ^{\prime
\prime }\left( 0\right) \right] ^{2}+\left[ \left( Lu\right) ^{(3)}\left(
0\right) \right] ^{2}\leq \left(
a_{1}^{2}+a_{2}^{2}+a_{3}^{2}+a_{4}^{2}\right) \left\Vert u\right\Vert _{W_{%
{\Large 2}}^{{\LARGE 5}}}^{2}.
\end{equation*}%
Since

\begin{equation*}
\left( Lu\right) ^{(4)}=\left\langle u,\left( LR_{x}\right)
^{(4)}\right\rangle _{W_{{\LARGE 2}}^{{\LARGE 5}}},
\end{equation*}%
then%
\begin{equation*}
\left\vert \left( Lu\right) ^{(4)}\right\vert \leq \left\Vert u\right\Vert
_{W_{2}^{5}}\left\Vert \left( LR_{x}\right) ^{(4)}\right\Vert
_{W_{2}^{5}}=a_{5}\left\Vert u\right\Vert _{W_{2}^{4}},\text{ \ \ (where }%
a_{5}>0\text{ is a positive constant),}
\end{equation*}%
so, we have

\begin{equation*}
\left[ \left( Lu\right) ^{(4)}\right] ^{2}\leq a_{5}^{2}\left\Vert
u\right\Vert _{{\LARGE W}_{{\LARGE 2}}^{{\LARGE 5}}}^{2},
\end{equation*}%
and

\begin{equation*}
\int_{0}^{1}\left[ \left( Lu\right) ^{(4)}\left( x\right) \right] ^{2}dx\leq
a_{5}^{2}\left\Vert u\right\Vert _{{\LARGE W}_{{\LARGE 2}}^{{\LARGE 5}}}^{2},
\end{equation*}%
that is%
\begin{eqnarray*}
\left\Vert Lu\right\Vert _{W_{2}^{4}}^{2} &=&\left[ \left( Lu\right) \left(
0\right) \right] ^{2}+\left[ \left( Lu\right) ^{\prime }\left( 0\right) %
\right] ^{2}+\left[ \left( Lu\right) ^{\prime \prime }\left( 0\right) \right]
^{2}+\left[ \left( Lu\right) ^{(3)}\left( 0\right) \right] ^{2}+\int_{0}^{1}%
\left[ \left( Lu\right) ^{(4)}\left( x\right) \right] ^{2}dx \\
&& \\
&\leq &\left( a_{1}^{2}+a_{2}^{2}+a_{3}^{2}+a_{4}^{2}+a_{5}^{2}\right)
\left\Vert u\right\Vert _{W_{2}^{5}}^{2}=P\left\Vert u\right\Vert
_{W_{2}^{4}}^{2},
\end{eqnarray*}%
where $P=\left( a_{1}^{2}+a_{2}^{2}+a_{3}^{2}+a_{4}^{2}+a_{5}^{2}\right) >0$
is a positive constant. This completes the proof. \ \ \ \ \ \ \ \ $\square $

\bigskip

\section{The Normal Orthogonal Function System of $W_{2}^{5}[a,b]$}

Let $\left\{ x_{i}\right\} _{i=1}^{\infty }$ be any dense set in $[0,1]$ and 
$\Psi _{x}(y)=L^{\ast }r_{x}(y),$ where $L^{\ast }$ is adjoint operator of $%
L $ and $r_{x}(y)$ is given by Remark 2.1. Furthermore, for simplicity let $%
\Psi _{i}(x)=\Psi _{x_{i}}(x),$ namely,

\begin{equation*}
\Psi _{i}(x)\overset{def}{=}\Psi _{x_{i}}(x)=L^{\ast }r_{x_{i}}(x).
\end{equation*}%
Now one can deduce following lemmas.

\bigskip

\textbf{Lemma 3.2. }$\left\{ \Psi _{i}(x)\right\} _{i=1}^{\infty }$ is
complete system of $W_{2}^{5}[0,1].$

\bigskip

\textbf{Proof: }For $u\in W_{2}^{5}[0,1]$, let $\left\langle u,\Psi
_{i}\right\rangle =0$ $(i=1,2,...),$ that is%
\begin{equation*}
\left\langle u,L^{\ast }r_{x_{i}}\right\rangle =(Lu)(x_{i})=0.
\end{equation*}%
Note that $\left\{ x_{i}\right\} _{i=1}^{\infty }$ is the dense set in $%
[0,1],$ therefore $(Lu)(x)=0.$ It follows that $u(x)=0$ from the existence
of $L^{-1}.$ This completes the proof. \ \ \ \ $\square $

\bigskip

\textbf{Lemma 3.3. }The following formula holds%
\begin{equation*}
\Psi _{i}(x)=\left( L\eta R_{x}(\eta )\right) \left( x_{i}\right) ,
\end{equation*}%
where the subscript $\eta $ of operator $L\eta $ indicates that the operator 
$L$ applies to function of $\eta .$

\bigskip

\textbf{Proof:}%
\begin{eqnarray*}
\Psi _{i}(x) &=&\left\langle \Psi _{i}(\xi ),R_{x}(\xi )\right\rangle
_{W_{2}^{5}[0,1]} \\
&=&\left\langle L^{\ast }r_{x_{i}}\left( \xi \right) ,R_{x}(\xi
)\right\rangle _{W_{2}^{5}[0,1]} \\
&=&\left\langle \left( r_{x_{i}}\right) \left( \xi \right) ,\left( L_{\eta
}R_{x}(\eta )\right) \left( \xi \right) \right\rangle _{W_{2}^{4}[0,1]} \\
&=&\left( L_{\eta }R_{x}(\eta )\right) \left( x_{i}\right) .
\end{eqnarray*}%
This completes the proof. \ \ \ \ \ \ \ \ $\square $

\bigskip

\textbf{Remark 3.1.} The orthonormal system $\left\{ \overline{\Psi }%
_{i}(x)\right\} _{i=1}^{\infty }$ of $W_{2}^{5}[0,1]$ can be derived from
Gram-Schmidt orthogonalization process of $\left\{ \Psi _{i}(x)\right\}
_{i=1}^{\infty },$%
\begin{equation}
\overline{\Psi }_{i}(x)=\sum_{k=1}^{i}\beta _{ik}\Psi _{k}(x),\text{ \ }%
(\beta _{ii}>0,\text{ \ }i=1,2,...)  \tag{3.7}
\end{equation}%
where $\beta _{ik}$ are orthogonal coefficients.

In the following, we will give the representation of the exact solution of
Eq.(1.17) in the reproducing kernel space $W_{2}^{5}[0,1].$

\bigskip

\subsection{\textbf{The Structure of the Solution and the Main Results}}

\bigskip

\textbf{Theorem 3.2.} If $u$ is the exact solution of (3.1) then

\begin{equation*}
u=\sum_{i=1}^{\infty }\sum_{k=1}^{i}\beta
_{ik}M(x_{k},u(x_{k}),u^{(3)}(x_{k}))\overline{\Psi }_{i}(x),
\end{equation*}%
where $\left\{ x_{i}\right\} _{i=1}^{\infty }$ is a dense set in $[0,1].$

\bigskip

\textbf{Proof: }From the (3.7) and uniqueness of solution of (3.1), we have%
\begin{eqnarray*}
u &=&\sum_{i=1}^{\infty }\left\langle u,\overline{\Psi }_{i}\right\rangle
_{W_{2}^{5}}\overline{\Psi }_{i}=\sum_{i=1}^{\infty }\sum_{k=1}^{i}\beta
_{ik}\left\langle u,L^{\ast }r_{x_{k}}\right\rangle _{W_{2}^{5}}\overline{%
\Psi }_{i} \\
&=&\sum_{i=1}^{\infty }\sum_{k=1}^{i}\beta _{ik}\left\langle
Lu,r_{x_{k}}\right\rangle _{W_{2}^{4}}\overline{\Psi }_{i}=\sum_{i=1}^{%
\infty }\sum_{k=1}^{i}\beta _{ik}\left\langle
M(x,u,u^{(3)}),r_{x_{k}}\right\rangle _{W_{2}^{4}}\overline{\Psi }_{i} \\
&=&\sum_{i=1}^{\infty }\sum_{k=1}^{i}\beta
_{ik}M(x_{k},u(x_{k}),u^{(3)}(x_{k}))\overline{\Psi }_{i}(x).
\end{eqnarray*}%
This completes the proof. \ \ \ \ \ \ \ \ \ \ \ \ \ \ \ \ \ \ $\square $

Now the approximate solution $u_{n}$ can be obtained by truncating the $n-$
term of the exact solution $u$ as%
\begin{equation*}
u_{n}=\sum_{i=1}^{n}\sum_{k=1}^{i}\beta _{ik}M(x_{k},u(x_{k}),u^{(3)}(x_{k}))%
\overline{\Psi }_{i}(x).
\end{equation*}

\bigskip

\textbf{Lemma 3.4. }$\left( \cite{ali1}\right) $Assume $u$ is the solution
of (3.1) and $r_{n}$ is the error between the approximate solution $u_{n}$
and the exact solution $u.$ Then the error sequence $r_{n}$ is monotone
decreasing in the sense of $\left\Vert .\right\Vert _{W_{2}^{5}}$ and $%
\left\Vert r_{n}(x)\right\Vert _{W_{2}^{5}}\rightarrow 0$.

\bigskip

\bigskip

\section{Numerical Results}

In this section, comparisons of results have been made through different
Reynolds numbers $\func{Re}$\ and magnetic field effect $m$.\ All
computations are performed by Maple 16. The RKM does not require
discretization of the variables, i.e., time and space, it is not effected by
computation round off errors and one is not faced with necessity of large
computer memory and time. The accuracy of the RKM for the MHD squeezing
fluid flow is controllable and absolute errors are small with present choice
of $x$ (see Tables 1-12). The numerical results we obtained justify the
advantage of this methodology.

\bigskip

\begin{center}
\begin{equation*}
\begin{tabular}{|l|l|l|l|l|l|}
\hline
$x$ & 
\begin{tabular}{l}
Numerical \\ 
Solution \\ 
$\left( RK-4\right) $%
\end{tabular}
& 
\begin{tabular}{l}
Approximate \\ 
Solution \\ 
\end{tabular}
& Absolute Error & Relative Error & Time (s) \\ \hline
$0.1$ & $0.150294$ & $0.15029400074386619072$ & $7.4386\times 10^{-10}$ & $%
4.9494071002\times 10^{-9}$ & $2.948$ \\ \hline
$0.2$ & $0.297481$ & $0.29748099943286204844$ & $5.6713\times 10^{-10}$ & $%
1.9064678132\times 10^{-9}$ & $2.980$ \\ \hline
$0.3$ & $0.438467$ & $0.43846699936146542481$ & $6.3853\times 10^{-10}$ & $%
1.4562887861\times 10^{-9}$ & $2.870$ \\ \hline
$0.4$ & $0.570189$ & $0.57018899983086605298$ & $1.6913\times 10^{-10}$ & $%
2.9662786728\times 10^{-10}$ & $2.792$ \\ \hline
$0.5$ & $0.689624$ & $0.68962399932753349664$ & $6.7246\times 10^{-10}$ & $%
9.7512050531\times 10^{-10}$ & $2.824$ \\ \hline
$0.6$ & $0.793796$ & $0.79379600052975674440$ & $5.2975\times 10^{-10}$ & $%
6.6737139569\times 10^{-10}$ & $2.902$ \\ \hline
$0.7$ & $0.879779$ & $0.87977900034152532706$ & $3.4152\times 10^{-10}$ & $%
3.8819445231\times 10^{-10}$ & $2.964$ \\ \hline
$0.8$ & $0.944696$ & $0.94469600021478585921$ & $2.1478\times 10^{-10}$ & $%
2.2735976357\times 10^{-10}$ & $2.808$ \\ \hline
$0.9$ & $0.985707$ & $0.98570699945336089741$ & $5.46639\times 10^{-10}$ & $%
5.5456550738\times 10^{-10}$ & $2.761$ \\ \hline
$1.0$ & $1.0$ & $1.0$ & $0.0$ & $0.0$ & $2.902$ \\ \hline
\end{tabular}%
\end{equation*}

\textbf{Table 4.1.} Numerical results at $m=1$ and $Re=1.$

\bigskip

\begin{equation*}
\begin{tabular}{|l|l|}
\hline
$x$ & 
\begin{tabular}{l}
Numerical \\ 
Solution \\ 
$\left( RK-4\right) $%
\end{tabular}
\\ \hline
$0.0$ & $0.0$ \\ \hline
$0.1$ & $0.150294$ \\ \hline
$0.2$ & $0.297481$ \\ \hline
$0.3$ & $0.438467$ \\ \hline
$0.4$ & $0.570189$ \\ \hline
$0.5$ & $0.689624$ \\ \hline
$0.6$ & $0.793796$ \\ \hline
$0.7$ & $0.879779$ \\ \hline
$0.8$ & $0.944696$ \\ \hline
$0.9$ & $0.985707$ \\ \hline
$1.0$ & $1.0$ \\ \hline
\end{tabular}%
\begin{tabular}{|l|}
\hline
\begin{tabular}{l}
\\ 
\\ 
OHAM%
\end{tabular}
\\ \hline
$0.0$ \\ \hline
$0.150265$ \\ \hline
$0.297424$ \\ \hline
$0.438387$ \\ \hline
$0.570093$ \\ \hline
$0.68952$ \\ \hline
$0.793695$ \\ \hline
$0.879695$ \\ \hline
$0.944641$ \\ \hline
$0.985687$ \\ \hline
$1.0$ \\ \hline
\end{tabular}%
\begin{tabular}{|l|}
\hline
\begin{tabular}{l}
\\ 
\\ 
RKHSM%
\end{tabular}
\\ \hline
$0.0$ \\ \hline
$0.15029400074386619072$ \\ \hline
$0.29748099943286204844$ \\ \hline
$0.43846699936146542481$ \\ \hline
$0.57018899983086605298$ \\ \hline
$0.68962399932753349664$ \\ \hline
$0.79379600052975674440$ \\ \hline
$0.87977900034152532706$ \\ \hline
$0.94469600021478585921$ \\ \hline
$0.98570699945336089741$ \\ \hline
$1.0$ \\ \hline
\end{tabular}%
\end{equation*}

\textbf{Table 4.2.} Comparison between RK-4, OHAM and RKHSM solutions at $%
m=1 $ and $Re=1.$

\bigskip
\end{center}

\begin{equation*}
\begin{tabular}{|l|l|l|l|l|l|}
\hline
$x$ & 
\begin{tabular}{l}
Numerical \\ 
Solution \\ 
$\left( RK-4\right) $%
\end{tabular}
& 
\begin{tabular}{l}
Approximate \\ 
Solution \\ 
\end{tabular}
& Absolute Error & Relative Error & Time (s) \\ \hline
$0.1$ & $0.137044$ & $0.13704399924397146430$ & $7.560285\times 10^{-10}$ & $%
5.51668468\times 10^{-9}$ & $3.261$ \\ \hline
$0.2$ & $0.272494$ & $0.27249400041809657591$ & $4.180965\times 10^{-10}$ & $%
1.53433314\times 10^{-9}$ & $3.542$ \\ \hline
$0.3$ & $0.404637$ & $0.40463699937791012358$ & $6.220898\times 10^{-10}$ & $%
1.53740235\times 10^{-9}$ & $2.949$ \\ \hline
$0.4$ & $0.531508$ & $0.53150799980699743080$ & $1.930025\times 10^{-10}$ & $%
3.63122604\times 10^{-10}$ & $3.541$ \\ \hline
$0.5$ & $0.650756$ & $0.65075599905912100256$ & $9.408789\times 10^{-10}$ & $%
1.44582454\times 10^{-9}$ & $3.089$ \\ \hline
$0.6$ & $0.759478$ & $0.75947799979255971384$ & $2.074402\times 10^{-10}$ & $%
2.73135345\times 10^{-10}$ & $2.996$ \\ \hline
$0.7$ & $0.854035$ & $0.85403499924057783299$ & $7.594221\times 10^{-10}$ & $%
8.89216679\times 10^{-10}$ & $3.026$ \\ \hline
$0.8$ & $0.929817$ & $0.92981700082221438640$ & $8.222143\times 10^{-10}$ & $%
8.84275493\times 10^{-10}$ & $7.582$ \\ \hline
$0.9$ & $0.980963$ & $0.98096299961587653980$ & $3.841234\times 10^{-10}$ & $%
3.91577929\times 10^{-10}$ & $3.291$ \\ \hline
$1.0$ & $1.0$ & $1.0$ & $0.0$ & $0.0$ & $2.902$ \\ \hline
\end{tabular}%
\end{equation*}

\begin{center}
\textbf{Table 4.3.} Numerical results at $m=3$ and $Re=1.$
\end{center}

\bigskip

\begin{equation*}
\begin{tabular}{|l|l|}
\hline
$x$ & 
\begin{tabular}{l}
Numerical \\ 
Solution \\ 
$\left( RK-4\right) $%
\end{tabular}
\\ \hline
$0.0$ & $0.0$ \\ \hline
$0.1$ & $0.137044$ \\ \hline
$0.2$ & $0.272494$ \\ \hline
$0.3$ & $0.404637$ \\ \hline
$0.4$ & $0.531508$ \\ \hline
$0.5$ & $0.650756$ \\ \hline
$0.6$ & $0.759478$ \\ \hline
$0.7$ & $0.854035$ \\ \hline
$0.8$ & $0.929817$ \\ \hline
$0.9$ & $0.980963$ \\ \hline
$1.0$ & $1.0$ \\ \hline
\end{tabular}%
\begin{tabular}{|l|}
\hline
\begin{tabular}{l}
\\ 
\\ 
OHAM%
\end{tabular}
\\ \hline
$0.0$ \\ \hline
$0.13709$ \\ \hline
$0.272583$ \\ \hline
$0.404759$ \\ \hline
$0.531649$ \\ \hline
$0.650894$ \\ \hline
$0.759591$ \\ \hline
$0.854106$ \\ \hline
$0.929845$ \\ \hline
$0.980966$ \\ \hline
$1.0$ \\ \hline
\end{tabular}%
\begin{tabular}{|l|}
\hline
\begin{tabular}{l}
\\ 
\\ 
RKHSM%
\end{tabular}
\\ \hline
$0.0$ \\ \hline
$0.13704399924397146430$ \\ \hline
$0.27249400041809657591$ \\ \hline
$0.40463699937791012358$ \\ \hline
$0.53150799980699743080$ \\ \hline
$0.65075599905912100256$ \\ \hline
$0.75947799979255971384$ \\ \hline
$0.85403499924057783299$ \\ \hline
$0.92981700082221438640$ \\ \hline
$0.98096299961587653980$ \\ \hline
$1.0$ \\ \hline
\end{tabular}%
\end{equation*}

\begin{center}
\textbf{Table 4.4.} Comparison between RK-4, OHAM and RKHSM solutions at $%
m=3 $ and $Re=1.$

\bigskip
\end{center}

\begin{equation*}
\begin{tabular}{|l|l|l|l|l|l|}
\hline
$x$ & 
\begin{tabular}{l}
Numerical \\ 
Solution \\ 
$\left( RK-4\right) $%
\end{tabular}
& 
\begin{tabular}{l}
Approximate \\ 
Solution \\ 
\end{tabular}
& Absolute Error & Relative Error & Time (s) \\ \hline
$0.1$ & $0.114976$ & $0.11497599095960418967$ & $9.040395\times 10^{-9}$ & $%
7.86285469\times 10^{-8}$ & $4.290$ \\ \hline
$0.2$ & $0.229882$ & $0.22988199268533318687$ & $7.314666\times 10^{-9}$ & $%
3.18192238\times 10^{-8}$ & $4.134$ \\ \hline
$0.3$ & $0.344604$ & $0.34460400584434350472$ & $5.844343\times 10^{-9}$ & $%
1.69595927\times 10^{-8}$ & $4.477$ \\ \hline
$0.4$ & $0.458904$ & $0.45890399132822355411$ & $8.671776\times 10^{-9}$ & $%
1.88967113\times 10^{-8}$ & $4.275$ \\ \hline
$0.5$ & $0.572276$ & $0.5722759999680104400$ & $3.198956\times 10^{-11}$ & $%
5.5898832\times 10^{-11}$ & $3.931$ \\ \hline
$0.6$ & $0.683628$ & $0.68362799155831029523$ & $8.441689\times 10^{-9}$ & $%
1.23483673\times 10^{-8}$ & $4.556$ \\ \hline
$0.7$ & $0.790607$ & $0.79060700783664672119$ & $7.836646\times 10^{-9}$ & $%
9.9121899\times 10^{-9}$ & $4.461$ \\ \hline
$0.8$ & $0.888173$ & $0.88817300466724146312$ & $4.667241\times 10^{-9}$ & $%
5.25487879\times 10^{-9}$ & $3.885$ \\ \hline
$0.9$ & $0.965578$ & $0.96557800220185786369$ & $2.201857\times 10^{-9}$ & $%
2.28035214\times 10^{-9}$ & $5.007$ \\ \hline
$1.0$ & $1.0$ & $1.0$ & $0.0$ & $0.0$ & $2.902$ \\ \hline
\end{tabular}%
\end{equation*}

\begin{center}
\textbf{Table 4.5.} Numerical results at $m=8$ and $Re=1.$
\end{center}

\bigskip

\begin{equation*}
\begin{tabular}{|l|l|}
\hline
$x$ & 
\begin{tabular}{l}
Numerical \\ 
Solution \\ 
$\left( RK-4\right) $%
\end{tabular}
\\ \hline
$0.0$ & $0.0$ \\ \hline
$0.1$ & $0.114976$ \\ \hline
$0.2$ & $0.229882$ \\ \hline
$0.3$ & $0.344604$ \\ \hline
$0.4$ & $0.458904$ \\ \hline
$0.5$ & $0.572276$ \\ \hline
$0.6$ & $0.683628$ \\ \hline
$0.7$ & $0.790607$ \\ \hline
$0.8$ & $0.888173$ \\ \hline
$0.9$ & $0.965578$ \\ \hline
$1.0$ & $1.0$ \\ \hline
\end{tabular}%
\begin{tabular}{|l|}
\hline
\begin{tabular}{l}
\\ 
\\ 
OHAM%
\end{tabular}
\\ \hline
$0.0$ \\ \hline
$0.11507$ \\ \hline
$0.230068$ \\ \hline
$0.344866$ \\ \hline
$0.459205$ \\ \hline
$0.572545$ \\ \hline
$0.683769$ \\ \hline
$0.790543$ \\ \hline
$0.887936$ \\ \hline
$0.965381$ \\ \hline
$1.0$ \\ \hline
\end{tabular}%
\begin{tabular}{|l|}
\hline
\begin{tabular}{l}
\\ 
\\ 
RKHSM%
\end{tabular}
\\ \hline
$0.0$ \\ \hline
$0.11497599095960418967$ \\ \hline
$0.22988199268533318687$ \\ \hline
$0.34460400584434350472$ \\ \hline
$0.45890399132822355411$ \\ \hline
$0.5722759999680104400$ \\ \hline
$0.68362799155831029523$ \\ \hline
$0.79060700783664672119$ \\ \hline
$0.88817300466724146312$ \\ \hline
$0.96557800220185786369$ \\ \hline
$1.0$ \\ \hline
\end{tabular}%
\end{equation*}

\begin{center}
\textbf{Table 4.6.} Comparison between RK-4, OHAM and RKHSM solutions at $%
m=8 $ and $Re=1.$

\bigskip
\end{center}

\begin{equation*}
\begin{tabular}{|l|l|l|l|l|l|}
\hline
$x$ & 
\begin{tabular}{l}
Numerical \\ 
Solution \\ 
$\left( RK-4\right) $%
\end{tabular}
& 
\begin{tabular}{l}
Approximate \\ 
Solution \\ 
\end{tabular}
& Absolute Error & Relative Error & Time (s) \\ \hline
$0.1$ & $0.105391$ & $0.10539098947593257979$ & $1.0524067\times 10^{-8}$ & $%
9.985736372\times 10^{-8}$ & $4.134$ \\ \hline
$0.2$ & $0.210782$ & $0.2107819933190829$ & $6.6809171\times 10^{-9}$ & $%
3.16958616\times 10^{-8}$ & $5.101$ \\ \hline
$0.3$ & $0.316173$ & $0.3161729190893567630$ & $8.0910643\times 10^{-8}$ & $%
2.559062387\times 10^{-7}$ & $3.010$ \\ \hline
$0.4$ & $0.421563$ & $0.4215629919618786430$ & $8.0381213\times 10^{-9}$ & $%
1.906742611\times 10^{-8}$ & $3.198$ \\ \hline
$0.5$ & $0.526952$ & $0.5269519479728988$ & $5.2027101\times 10^{-8}$ & $%
9.873214486\times 10^{-8}$ & $3.042$ \\ \hline
$0.6$ & $0.632324$ & $0.632323981769674315$ & $1.8230325\times 10^{-8}$ & $%
2.883067175\times 10^{-8}$ & $3.074$ \\ \hline
$0.7$ & $0.737586$ & $0.7375860570172070642$ & $5.7017207\times 10^{-8}$ & $%
7.730245295\times 10^{-8}$ & $3.089$ \\ \hline
$0.8$ & $0.842051$ & $0.84205103495023398982$ & $3.4950233\times 10^{-8}$ & $%
4.150607741\times 10^{-8}$ & $3.073$ \\ \hline
$0.9$ & $0.940861$ & $0.94086101815219431313$ & $1.8152194\times 10^{-8}$ & $%
1.929317328\times 10^{-8}$ & $3.135$ \\ \hline
$1.0$ & $1.0$ & $1.0$ & $0.0$ & $0.0$ & $2.902$ \\ \hline
\end{tabular}%
\end{equation*}

\begin{center}
\textbf{Table 4.7.} Numerical results at $m=20$ and $Re=1.$
\end{center}

\bigskip

\begin{equation*}
\begin{tabular}{|l|l|}
\hline
$x$ & 
\begin{tabular}{l}
Numerical \\ 
Solution \\ 
$\left( RK-4\right) $%
\end{tabular}
\\ \hline
$0.0$ & $0.0$ \\ \hline
$0.1$ & $0.105391$ \\ \hline
$0.2$ & $0.210782$ \\ \hline
$0.3$ & $0.316173$ \\ \hline
$0.4$ & $0.421563$ \\ \hline
$0.5$ & $0.526952$ \\ \hline
$0.6$ & $0.632324$ \\ \hline
$0.7$ & $0.737586$ \\ \hline
$0.8$ & $0.842051$ \\ \hline
$0.9$ & $0.940861$ \\ \hline
$1.0$ & $1.0$ \\ \hline
\end{tabular}%
\begin{tabular}{|l|}
\hline
\begin{tabular}{l}
\\ 
\\ 
OHAM%
\end{tabular}
\\ \hline
$0.0$ \\ \hline
$0.105312$ \\ \hline
$0.210625$ \\ \hline
$0.315938$ \\ \hline
$0.421249$ \\ \hline
$0.526551$ \\ \hline
$0.631824$ \\ \hline
$0.736971$ \\ \hline
$0.841352$ \\ \hline
$0.94035$ \\ \hline
$1.0$ \\ \hline
\end{tabular}%
\begin{tabular}{|l|}
\hline
\begin{tabular}{l}
\\ 
\\ 
RKHSM%
\end{tabular}
\\ \hline
$0.0$ \\ \hline
$0.10539098947593257979$ \\ \hline
$0.2107819933190829$ \\ \hline
$0.3161729190893567630$ \\ \hline
$0.4215629919618786430$ \\ \hline
$0.5269519479728988$ \\ \hline
$0.632323981769674315$ \\ \hline
$0.7375860570172070642$ \\ \hline
$0.84205103495023398982$ \\ \hline
$0.94086101815219431313$ \\ \hline
$1.0$ \\ \hline
\end{tabular}%
\end{equation*}

\begin{center}
\textbf{Table 4.8.} Comparison between RK-4, OHAM and RKHSM solutions at $%
m=20$ and $Re=1.$

\bigskip
\end{center}

\begin{equation*}
\begin{tabular}{|l|l|l|l|l|l|}
\hline
$x$ & 
\begin{tabular}{l}
Numerical \\ 
Solution \\ 
$\left( RK-4\right) $%
\end{tabular}
& 
\begin{tabular}{l}
Approximate \\ 
Solution \\ 
\end{tabular}
& Absolute Error & Relative Error & Time (s) \\ \hline
$0.1$ & $0.158104$ & $0.15810400012535311729$ & $1.25353117\times 10^{-10}$
& $7.928522826\times 10^{-10}$ & $5.304$ \\ \hline
$0.2$ & $0.311962$ & $0.31196200057873017887$ & $5.78730178\times 10^{-10}$
& $1.855130364\times 10^{-9}$ & $7.332$ \\ \hline
$0.3$ & $0.457539$ & $0.45753900003164153289$ & $3.16415328\times 10^{-11}$
& $6.915592526\times 10^{-11}$ & $5.913$ \\ \hline
$0.4$ & $0.591193$ & $0.59119300033029000468$ & $3.30290004\times 10^{-10}$
& $5.586838894\times 10^{-10}$ & $6.272$ \\ \hline
$0.5$ & $0.709771$ & $0.70977100026331200670$ & $2.63312006\times 10^{-10}$
& $3.709816359\times 10^{-10}$ & $5.757$ \\ \hline
$0.6$ & $0.810642$ & $0.81064200064720692438$ & $6.47206924\times 10^{-10}$
& $7.983880978\times 10^{-10}$ & $6.256$ \\ \hline
$0.7$ & $0.891666$ & $0.89166599939606220359$ & $6.03937796\times 10^{-10}$
& $6.039377964\times 10^{-10}$ & $6.396$ \\ \hline
$0.8$ & $0.95112$ & $0.95112000044608660232$ & $4.46086602\times 10^{-10}$ & 
$4.690119035\times 10^{-10}$ & $5.101$ \\ \hline
$0.9$ & $0.987612$ & $0.98761199979328069240$ & $2.06719307\times 10^{-10}$
& $2.093122679\times 10^{-10}$ & $5.616$ \\ \hline
$1.0$ & $1.0$ & $1.0$ & $0.0$ & $0.0$ & $2.902$ \\ \hline
\end{tabular}%
\end{equation*}

\begin{center}
\textbf{Table 4.9.} Numerical results at $m=1$ and $Re=4.$
\end{center}

\bigskip

\begin{equation*}
\begin{tabular}{|l|l|}
\hline
$x$ & 
\begin{tabular}{l}
Numerical \\ 
Solution \\ 
$\left( RK-4\right) $%
\end{tabular}
\\ \hline
$0.0$ & $0.0$ \\ \hline
$0.1$ & $0.158104$ \\ \hline
$0.2$ & $0.311962$ \\ \hline
$0.3$ & $0.457539$ \\ \hline
$0.4$ & $0.591193$ \\ \hline
$0.5$ & $0.709771$ \\ \hline
$0.6$ & $0.810642$ \\ \hline
$0.7$ & $0.891666$ \\ \hline
$0.8$ & $0.95112$ \\ \hline
$0.9$ & $0.987612$ \\ \hline
$1.0$ & $1.0$ \\ \hline
\end{tabular}%
\begin{tabular}{|l|}
\hline
\\ \hline
\\ \hline
OHAM \\ \hline
$0.0$ \\ \hline
$0.156218$ \\ \hline
$0.308363$ \\ \hline
$0.452557$ \\ \hline
$0.585287$ \\ \hline
$0.703518$ \\ \hline
$0.804726$ \\ \hline
$0.886838$ \\ \hline
$0.948051$ \\ \hline
$0.986529$ \\ \hline
$1.0$ \\ \hline
\end{tabular}%
\begin{tabular}{|l|}
\hline
\\ \hline
\\ \hline
RKHSM \\ \hline
$0.0$ \\ \hline
$0.15810400012535311729$ \\ \hline
$0.31196200057873017887$ \\ \hline
$0.45753900003164153289$ \\ \hline
$0.59119300033029000468$ \\ \hline
$0.70977100026331200670$ \\ \hline
$0.81064200064720692438$ \\ \hline
$0.89166599939606220359$ \\ \hline
$0.95112000044608660232$ \\ \hline
$0.98761199979328069240$ \\ \hline
$1.0$ \\ \hline
\end{tabular}%
\end{equation*}

\begin{center}
\textbf{Table 4.10.} Comparison between RK-4, OHAM and RKHSM solutions at $%
m=1$ and $Re=4.$

\bigskip 
\end{center}

\begin{equation*}
\begin{tabular}{|l|l|l|l|l|l|}
\hline
$x$ & 
\begin{tabular}{l}
Numerical \\ 
Solution \\ 
$\left( RK-4\right) $%
\end{tabular}
& 
\begin{tabular}{l}
Approximate \\ 
Solution \\ 
\end{tabular}
& Absolute Error & Relative Error & Time (s) \\ \hline
$0.1$ & $0.167616$ & $0.1676160001397322991$ & $1.39732299\times 10^{-10}$ & 
$8.3364535\times 10^{-10}$ & $5.569$ \\ \hline
$0.2$ & $0.329031$ & $0.32903100221406728329$ & $2.21406728\times 10^{-9}$ & 
$6.7290537\times 10^{-9}$ & $6.365$ \\ \hline
$0.3$ & $0.478907$ & $0.47890699791462877619$ & $2.08537122\times 10^{-9}$ & 
$4.3544388\times 10^{-9}$ & $7.378$ \\ \hline
$0.4$ & $0.613252$ & $0.61325199550552162812$ & $4.49447837\times 10^{-9}$ & 
$7.3289257\times 10^{-9}$ & $7.254$ \\ \hline
$0.5$ & $0.729428$ & $0.72942799845508679063$ & $1.5449132\times 10^{-9}$ & $%
2.117979\times 10^{-9}$ & $6.271$ \\ \hline
$0.6$ & $0.825843$ & $0.82584300690485584332$ & $6.9048558\times 10^{-9}$ & $%
8.3609788\times 10^{-9}$ & $7.425$ \\ \hline
$0.7$ & $0.901576$ & $0.90157600840425340903$ & $8.4042534\times 10^{-9}$ & $%
9.32173594\times 10^{-9}$ & $6.162$ \\ \hline
$0.8$ & $0.901576$ & $0.90157518382496567601$ & $8.16175\times 10^{-7}$ & $%
9.052759\times 10^{-7}$ & $7.410$ \\ \hline
$0.9$ & $0.988978$ & $0.98897799997420425356$ & $2.579574\times 10^{-11}$ & $%
2.6083235\times 10^{-11}$ & $7.910$ \\ \hline
$1.0$ & $1.0$ & $1.0$ & $0.0$ & $0.0$ & $2.902$ \\ \hline
\end{tabular}%
\end{equation*}

\begin{center}
\textbf{Table 4.11.} Numerical results at $m=1$ and $Re=10.$
\end{center}

\bigskip

\begin{equation*}
\begin{tabular}{|l|l|}
\hline
$x$ & 
\begin{tabular}{l}
Numerical \\ 
Solution \\ 
$\left( RK-4\right) $%
\end{tabular}
\\ \hline
$0.0$ & $0.0$ \\ \hline
$0.1$ & $0.167616$ \\ \hline
$0.2$ & $0.329031$ \\ \hline
$0.3$ & $0.478907$ \\ \hline
$0.4$ & $0.613252$ \\ \hline
$0.5$ & $0.729428$ \\ \hline
$0.6$ & $0.825843$ \\ \hline
$0.7$ & $0.901576$ \\ \hline
$0.8$ & $0.901576$ \\ \hline
$0.9$ & $0.988978$ \\ \hline
$1.0$ & $1.0$ \\ \hline
\end{tabular}%
\begin{tabular}{|l|}
\hline
\\ \hline
\\ \hline
OHAM \\ \hline
$0.0$ \\ \hline
$0.175911$ \\ \hline
$0.344336$ \\ \hline
$0.498671$ \\ \hline
$0.633941$ \\ \hline
$0.747277$ \\ \hline
$0.838004$ \\ \hline
$0.907244$ \\ \hline
$0.956954$ \\ \hline
$0.988387$ \\ \hline
$1.0$ \\ \hline
\end{tabular}%
\begin{tabular}{|l|}
\hline
\\ \hline
\\ \hline
RKHSM \\ \hline
$0.0$ \\ \hline
$0.1676160001397322991$ \\ \hline
$0.32903100221406728329$ \\ \hline
$0.47890699791462877619$ \\ \hline
$0.61325199550552162812$ \\ \hline
$0.72942799845508679063$ \\ \hline
$0.82584300690485584332$ \\ \hline
$0.90157600840425340903$ \\ \hline
$0.90157518382496567601$ \\ \hline
$0.98897799997420425356$ \\ \hline
$1.0$ \\ \hline
\end{tabular}%
\end{equation*}

\begin{center}
\textbf{Table 4.12.} Comparison between RK-4, OHAM and RKHSM solutions at $%
m=1$ and $Re=10.$
\end{center}

\bigskip

\section{Conclusion}

In this paper, we introduced an algorithm for solving the MHD squeezing
fluid flow with boundary conditions. The method gives more realistic series
solutions that converge very rapidly in physical problems. The approximate
solution obtained by the present method is uniformly convergent.\newline

Clearly, the series solution methodology can be applied to much more
complicated nonlinear differential equations and boundary value problems.
However, if the problem becomes nonlinear, then the RKM does not require
discretization or perturbation and it does not make closure approximation.
Results of numerical examples show that the present method is an accurate
and reliable analytical method for the MHD squeezing fluid flow problem with
boundary conditions.

\begin{center}
$\bigskip $

$\bigskip $
\end{center}

Mustafa Inc: Department of Mathematics, Science Faculty, F\i rat University,
23119 Elaz\i \u{g} / Turkey, minc@firat.edu.tr

Ali Akg\"{u}l: Department of Mathematics, Education Faculty, Dicle
University, 21280 Diyarbak\i r / Turkey, aliakgul00727@gmail.com

\end{document}